\documentclass[11pt,a4paper]{amsart}
\usepackage[margin=30mm]{geometry}
\usepackage{amsmath,amssymb,amsthm,mathtools}
\usepackage{enumitem}
\usepackage{xcolor}
\usepackage{hyperref}
\usepackage{tikz}
\usepackage{float}
\hypersetup{hidelinks}
\newtheorem{theorem}{Theorem}[section]
\newtheorem{proposition}[theorem]{Proposition}
\newtheorem{lemma}[theorem]{Lemma}
\newtheorem{corollary}[theorem]{Corollary}
\theoremstyle{definition}

\newtheorem{example}[theorem]{Example}
\theoremstyle{remark}

\newcommand{\ZZ}{\mathbb Z}
\newcommand{\RR}{\mathbb R}
\newcommand{\calA}{\mathcal A}

\newcommand{\calC}{\mathcal C}

\newcommand{\sgn}{\operatorname{sgn}}

\numberwithin{equation}{section}

\title[Chromatic polynomials of signed graphs]
{Chromatic polynomials of signed graphs with a balancing vertex}

\author{Akihiro Higashitani}
\address{Department of Pure and Applied Mathematics, Graduate School of Information Science and Technology, Osaka University, Suita, Osaka 565-0871, Japan}
\email{higashitani@ist.osaka-u.ac.jp}

\author{Go Kamata}
\address{Department of Pure and Applied Mathematics, Graduate School of Information Science and Technology, Osaka University, Suita, Osaka 565-0871, Japan}
\email{kamata.gou@ist.osaka-u.ac.jp}

\author{Masato Konoike}
\address{Department of Pure and Applied Mathematics, Graduate School of Information Science and Technology, Osaka University, Suita, Osaka 565-0871, Japan}
\email{kounoike-m@ist.osaka-u.ac.jp}


\begin{document}

\begin{abstract}
We study the zero-free chromatic polynomial $Z_\Sigma(q)$ of a signed graph $\Sigma$ and the characteristic polynomial $\chi_\Sigma(q)=Z_\Sigma(q-1)$ of its associated hyperplane arrangement containing all coordinate hyperplanes. For a connected signed graph on $n\ge 2$ vertices with a balancing vertex, we prove the equivalence of the following four conditions: 
(i) $Z_\Sigma(q)$ is the chromatic polynomial of an ordinary graph; 
(ii) the numerator of the generating function of of the zero-free chromatic polynomial associated with $\Sigma$ has nonnegative coefficients; 
(iii) $\Sigma$ has a balanced block; and 
(iv) $\chi_\Sigma(2)=0$. We also provide explicit counterexamples showing that neither the balancing-vertex assumption nor connectedness can be omitted.
\end{abstract}

\maketitle

\section{Introduction}

The chromatic polynomial $P_G(q)$ of an ordinary graph $G$ counts its proper colorings with $q$ colors.
Signed graphs provide a natural extension of this theory: a positive edge requires distinct colors at its endpoints, whereas a negative edge requires colors that are not opposites.
Zaslavsky developed the corresponding coloring theory, including chromatic polynomials, switching invariance,
and their relation to hyperplane arrangements and signed-graphic matroids \cite{ZaslavskyColoring,ZaslavskyInvariants}.
An essential feature is the distinction between colorings with the color set $\{0,\pm1,\ldots,\pm k\}$ and zero-free colorings with the color set $\{\pm1,\ldots,\pm k\}$.
These give different polynomials in general. 
Chromatic invariants of signed graphs have also been investigated recently; see, e.g., \cite{GreavesSyatriadiUtomo}. 

Throughout this paper, let $\Sigma=(V,E^+,E^-)$ denote a loopless signed graph on the vertex set $V$ with the set of positive edges $E^+$ and negative edges $E^-$. 
We associate the following hyperplane arrangement in $\RR^V$ to $\Sigma$:
\begin{equation*}
 \calA_\Sigma = \{x_i=0:i\in V\} \cup\{x_i-x_j=0:ij\in E^+\} \cup\{x_i+x_j=0:ij\in E^-\}.
\end{equation*}
We use the following notation throughout the paper:
\begin{itemize}
\item $P_G(q)$ denotes the chromatic polynomial of an ordinary graph $G$, counting its proper colorings with $q$ colors.
\item $Z_\Sigma(q)$ denotes the zero-free chromatic polynomial of a signed graph $\Sigma$, where $Z_\Sigma(2k)$ counts proper zero-free colorings with colors in $\{\pm1,\ldots,\pm k\}$.
\item $\chi_\Sigma(q)$ denotes the characteristic polynomial of the hyperplane arrangement $\calA_\Sigma$.
For an ordinary graph $G$, viewed as a signed graph having all positive edges, we write $\chi_G(q)$ for the characteristic polynomial of $\calA_G$.
\end{itemize}
These polynomials satisfy
\begin{equation}\label{eq:shift}
 \chi_\Sigma(q)=Z_\Sigma(q-1).
\end{equation}
In particular, for an ordinary graph $G$, we have \[ \chi_G(q)=P_G(q-1).\]
Note that the shift comes from our convention that $\calA_\Sigma$ contains all coordinate hyperplanes.


For an ordinary graph $G$ on $n$ vertices, consider the polynomial $W_G(t)$ defined by
\[ \sum_{q\ge0}P_G(q)t^q=\frac{W_G(t)}{(1-t)^{n+1}}.\]
It is known that $W_G(t)$ has nonnegative integer coefficients.
This nonnegativity is closely related to Stanley's work on graph colorings and acyclic orientations \cite{Stanley}, and was established explicitly by Linial \cite{Linial}.
Steingr\'imsson \cite{Steingrimsson} introduced the coloring ideal and the coloring complex of a graph, connecting graph colorings with Stanley--Reisner theory.
In particular, $W_G(t)$ admits the following combinatorial interpretation:
\[ W_G(t)=\sum_{\pi\in S_n}t^{c_G(\pi)},\]
where $c_G(\pi)$ denotes the number of cuts of a permutation $\pi$ associated with $G$; see \cite[Theorem~6]{Steingrimsson}.
Further studies of chromatic generating functions include the work of Gansner and Vo \cite{GansnerVo}.
More recently, Beck and Le\'on \cite{BeckLeon} established inequalities 
for the coefficients in the corresponding binomial basis and investigated necessary conditions for a polynomial to be a chromatic polynomial.

Motivated by these results, we study the analogous generating function for signed graphs.
For a signed graph $\Sigma$ on $n$ vertices, we define $h_\Sigma(t)\in\ZZ[t]$ by
\begin{equation}\label{eq:series}
 \sum_{q\ge1}\chi_\Sigma(q)t^q
 =
 \frac{h_\Sigma(t)}{(1-t)^{n+1}}.
\end{equation}
If $\chi_\Sigma(q)=P_G(q-1)$ for some ordinary graph $G$, then $h_\Sigma(t)=tW_G(t)$, so $h_\Sigma(t)$ has nonnegative coefficients.
However, for general signed graphs, the converse does not hold. 
The main purpose of this paper is to characterize a class of signed graphs for which this nonnegativity condition is also sufficient 
for their characteristic polynomials to arise from ordinary graph colorings.

To state our main result, we recall some terminology concerning signed graphs.
A signed graph is \emph{balanced} if every cycle contains an even number of negative edges in its underlying graph. 
A \emph{block} of $\Sigma$ means a block of its underlying graph, equipped with the inherited edge signs. 
We regard every bridge as a block, and hence as a balanced block.
We call a vertex $u$ a \emph{balancing vertex} if $\Sigma \setminus u$ is balanced; this notion goes back to Zaslavsky \cite{ZaslavskyLocalized}. 
The existence of a balancing vertex in an unbalanced signed graph is equivalent to the condition that, after switching, 
all negative edges are incident to a common vertex; see \cite[Proposition~1.3]{ZaslavskyNegativeCircles}. 
See Subsection~\ref{sec:switching} for the details of switching. 

Splitting a balancing vertex is a standard construction in the theory of signed and biased graphs; see \cite[Section~2.1]{DeVosFunk}.
Our main result characterizes the nonnegativity of $h_\Sigma(t)$ and its relation to ordinary graph colorings in terms of the block structure of $\Sigma$.
The following is the main theorem of this paper. 

\begin{theorem}\label{thm:main}
Let $\Sigma$ be a connected loopless signed graph on $n \geq 2$ vertices with a balancing vertex. 
Then the following conditions are equivalent:
\begin{enumerate}[label=\textup{(\roman*)}]
 \item\label{main:graphic} There exists an ordinary graph $G$ such that
 \[ \chi_\Sigma(q)=\chi_G(q),\qquad\text{equivalently,}\quad Z_\Sigma(q)=P_G(q). \]
 \item\label{main:positive} Every coefficient of $h_\Sigma(t)$ is  nonnegative.
 \item\label{main:block} $\Sigma$ has a balanced block.
 \item\label{main:uniform} $Z_\Sigma(1) = 0$, equivalently, $\chi_\Sigma(2)=0$.
\end{enumerate}
Moreover, whenever these conditions hold, the graph $G$ can be chosen to be connected, simple, and to have $|V|$ vertices.
\end{theorem}

Recall from \cite[Corollary~3.3]{ZaslavskySignedGraphs} that a signed graph is balanced if and only if it is switching-equivalent to a signed graph having all positive edges. 
In particular, when $\Sigma$ is balanced, all four conditions (i)--(iv) of Theorem~\ref{thm:main} are satisfied. 
Thus, the main difficulty lies in the unbalanced case.
The implication \ref{main:graphic}$\Rightarrow$\ref{main:positive} follows directly from the aforementioned nonnegativity result for ordinary graphs \cite[Theorem~6]{Steingrimsson}.
The significance of Theorem~\ref{thm:main} is that this necessary condition becomes sufficient under the assumption of having a balancing vertex. 
Moreover, Theorem~\ref{thm:main} provides a characterization in terms of balanced blocks.

We briefly explain the ideas behind the proof of Theorem~\ref{thm:main}.
The implication \ref{main:positive}$\Rightarrow$\ref{main:block} is established by examining the first few coefficients of $h_\Sigma(t)$.
In fact, if $\Sigma$ has no balanced block, then the coefficient of $t^2$ or that of $t^3$ is negative, according to the parity of $n$.
This follows from positivity of the beta invariant \cite{Crapo,GreeneZaslavsky}, applied to ordinary graphs obtained by splitting the unbalanced blocks.
To prove \ref{main:block}$\Rightarrow$\ref{main:graphic}, we use a balanced block to construct an ordinary graph whose chromatic polynomial coincides with $Z_\Sigma(q)$.
The equivalence of \ref{main:block} and \ref{main:uniform} follows from \eqref{eq:block-product} and the fact that $Z_\Sigma(1)\ne0$ whenever $\Sigma$ has no balanced block. 
Finally, we give examples showing that the assumptions in Theorem~\ref{thm:main}, particularly the existence of a balancing vertex, cannot be omitted.

The paper is organized as follows. Section~\ref{sec:preliminaries} introduces the coloring conventions, establishes the required product and positivity formulas, 
and treats vertex splitting and its relation to the beta invariant. 
The main theorem is proved in Section~\ref{sec:proof}. 
Section~\ref{sec:examples} contains examples and counterexamples showing the role of the assumptions of Theorem~\ref{thm:main}. 

\bigskip

\section{Signed colorings, generating functions and vertex splitting}\label{sec:preliminaries}

In this section, we introduce the necessary conventions and establish lemmas for the proof of the main theorem.

We allow one positive edge and one negative edge between the same pair of vertices, but no repeated edges of the same sign.
Thus, $E^+$ and $E^-$ are subsets of $\binom{V}{2}$ and need not be disjoint.
The underlying multigraph retains both edges when both signs occur between the same pair of vertices, and cycles of length two are included in the definition of balance.
An \emph{ordinary graph} will always mean a finite simple graph. Thus, in particular, the main result includes signed graphs with simple underlying graphs.
A \emph{block} is a maximal connected subgraph with at least one edge and no cut vertex.
Every bridge forms a block. A two-vertex block may also consist of a pair of opposite signed parallel edges.
A graph with at least three vertices is \emph{$2$-connected} if it is connected and remains connected after deleting any one vertex.

\subsection{Zero-free colorings and switching}\label{sec:switching}

For $k\ge1$, let $\calC_k=\{\pm1,\ldots,\pm k\}$. A map $c:V\longrightarrow\calC_k$ is a \emph{proper zero-free coloring} of $\Sigma$ if
\[ c(i)\ne c(j)\quad(ij\in E^+), \qquad c(i)\ne-c(j)\quad(ij\in E^-).\]
The \emph{zero-free chromatic polynomial} $Z_\Sigma(q)$ of $\Sigma$ is defined by the condition
that $Z_\Sigma(2k)$ counts the number of proper zero-free colorings of $\Sigma$ with colors in $\calC_k$ for every $k\ge1$.
It is well known that $Z_\Sigma(q)$ is a monic polynomial in $\ZZ[q]$ of degree $n$ and satisfies $Z_\Sigma(0)=0$ for $n\ge1$; see \cite{ZaslavskyColoring,ZaslavskyInvariants}.
The identity \eqref{eq:shift} follows from the finite field method for hyperplane arrangements \cite{Athanasiadis};
see also \cite{ZaslavskyColoring,ZaslavskyInvariants} for the signed graph coloring interpretation.

A \emph{switching} is a map $\varepsilon:V\to\{1,-1\}$. 
Switching by $\varepsilon$ changes the sign $\sgn(e)$ of each edge $e=ij$ according to $\sgn(e)\mapsto\varepsilon(i)\sgn(e)\varepsilon(j)$.
The map $c(i)\mapsto\varepsilon(i)c(i)$ gives a bijection of proper zero-free colorings before and after switching.
Thus, $Z_\Sigma$, $\chi_\Sigma$, and $h_\Sigma$ are invariant under switching.
The product of the signs around a cycle is unchanged as well.

It is known that a signed graph is balanced if and only if it is switching-equivalent to a signed graph having all positive edges \cite[Corollary~3.3]{ZaslavskySignedGraphs}. 
Since switching preserves the zero-free chromatic polynomial, if $B$ is a balanced signed graph and $\overline B$ is its underlying ordinary graph, then
\begin{equation}\label{eq:balanced-ordinary}
 Z_B(q)=P_{\overline B}(q).
\end{equation}
Moreover, applying the same characterization to $\Sigma \setminus u$,
we see that $\Sigma \setminus u$ is balanced if and only if $\Sigma$ can be transformed by switching so that every negative edge is incident to $u$.

\subsection{Products over blocks}

\begin{lemma}\label{lem:vertex-sum}
The following statements hold.
\begin{enumerate}[label=\textup{(\alph*)}]
\item Suppose that signed graphs $\Sigma_1$ and $\Sigma_2$ have exactly one common vertex and no common edges.
Then $\displaystyle Z_{\Sigma_1\cup\Sigma_2}(q) =\frac{Z_{\Sigma_1}(q)Z_{\Sigma_2}(q)}{q}$.
Consequently, if a connected signed graph $\Sigma$ has blocks $B_1,\ldots,B_b$, then
\begin{equation}\label{eq:block-product}
 Z_\Sigma(q)=q^{1-b}\prod_{i=1}^b Z_{B_i}(q).
\end{equation}
\item Suppose that ordinary graphs $G_1$ and $G_2$ meet exactly in an edge and its two endpoints. Then
\begin{equation}\label{eq:edge-sum}
 P_{G_1\cup G_2}(q) =\frac{P_{G_1}(q)P_{G_2}(q)}{q(q-1)}.
\end{equation}
\end{enumerate}
\end{lemma}
\begin{proof}
(a) The permutations of $\calC_k$ commuting with color negation act transitively on $\calC_k$ and preserve proper zero-free colorings.
Hence, a prescribed color at the common vertex extends to $Z_{\Sigma_i}(2k)/(2k)$ proper colorings of $\Sigma_i$ for $i=1,2$.
Since the extensions on the two sides are independent, summing over the $2k$ possible common colors proves the first identity in (a) as a polynomial identity.
Iterating this identity along the block decomposition yields \eqref{eq:block-product}.

\noindent (b) For an integer $q\ge2$, every ordered pair of distinct colors on the common edge has the same number of extensions.
Since there are $q(q-1)$ such pairs, the formula \eqref{eq:edge-sum} follows. 
\end{proof}

\subsection{The first numerator coefficients}

Since $Z_\Sigma(0)=0$, \eqref{eq:shift} and \eqref{eq:series} imply that
\[ \sum_{q\ge1}\chi_\Sigma(q)t^q =t\sum_{q\ge0}Z_\Sigma(q)t^q.\]
Thus, $\deg h_\Sigma\le n+1$, and we may write $h_\Sigma(t)=\sum_{i=1}^{n+1}h_i t^i$.
For $n\ge2$, by multiplying \eqref{eq:series} by $(1-t)^{n+1}$, we obtain that 
\begin{equation}\label{eq:first-coefficients}
 h_1=0,\qquad h_2=\chi_\Sigma(2)=Z_\Sigma(1),\qquad h_3=Z_\Sigma(2)-(n+1)Z_\Sigma(1),
\end{equation}
where $Z_\Sigma(1)$ is a polynomial evaluation, not a count of colorings from a zero-free signed color set of odd cardinality.

\begin{lemma}\label{lem:two-colors}
If $\Sigma$ is connected, then $Z_\Sigma(2)\in\{0,2\}$.
\end{lemma}
\begin{proof}
With colors $\pm 1$, a positive edge forces opposite colors and a negative edge forces equal colors.
After the color at one vertex is fixed, the color at every vertex is forced by connectedness.
The requirements are either inconsistent, or give two colorings related by negating all colors.
\end{proof}

\subsection{Splitting a balancing vertex}\label{sec:splitting}

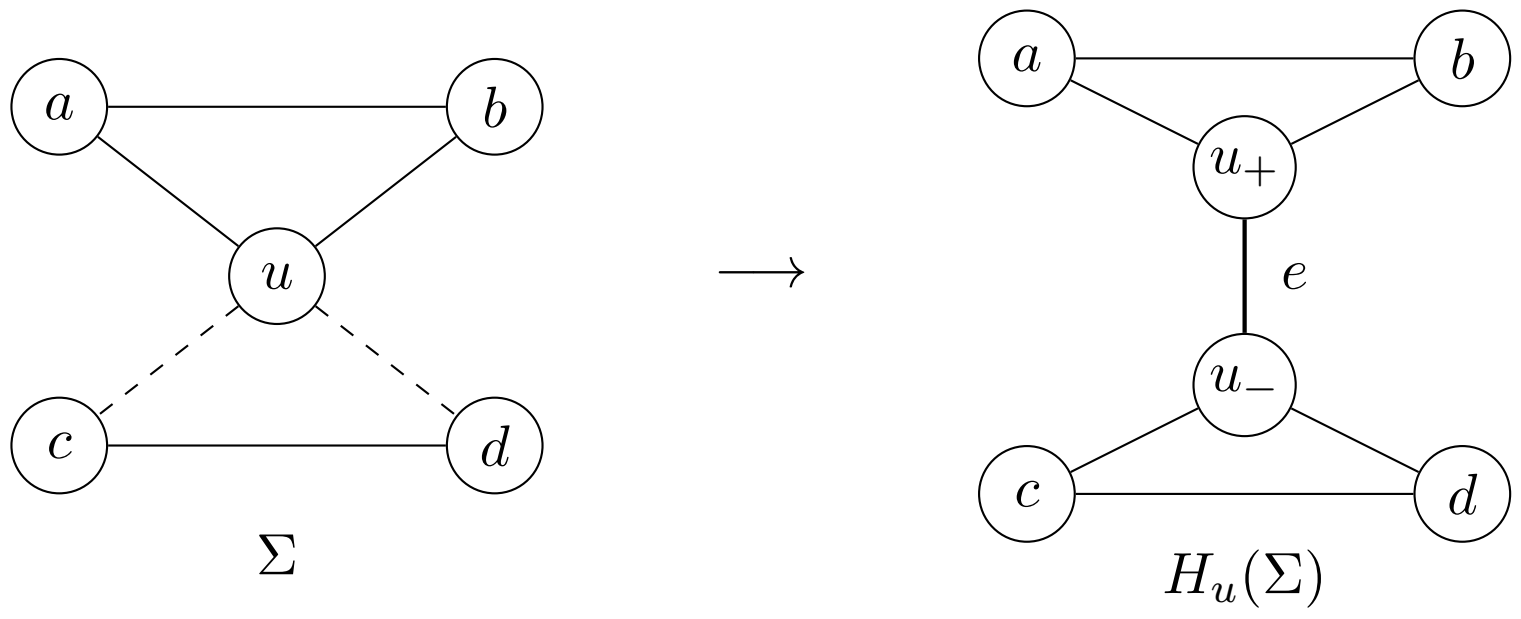
\begin{figure}[H]
\centering
\begin{tikzpicture}[scale=0.8, every node/.style={circle, draw, inner sep=1.5pt, minimum size=18pt}]

\begin{scope}[xshift=0cm]
\node (u) at (0,0) {$u$};
\node (a) at (-1.8,1.4) {$a$};
\node (b) at (1.8,1.4) {$b$};
\node (c) at (-1.8,-1.4) {$c$};
\node (d) at (1.8,-1.4) {$d$};

\draw (u) -- (a);
\draw (u) -- (b);
\draw (a) -- (b);
\draw (c) -- (d);

\draw[dashed] (u) -- (c);
\draw[dashed] (u) -- (d);

\node[draw=none, shape=rectangle] at (0,-2.3) {$\Sigma$};
\end{scope}

\node[draw=none, shape=rectangle] at (4,0) {$\longrightarrow$};

\begin{scope}[xshift=8cm]
\node (up) at (0,0.9) {$u_+$};
\node (um) at (0,-0.9) {$u_-$};
\node (a2) at (-1.8,1.8) {$a$};
\node (b2) at (1.8,1.8) {$b$};
\node (c2) at (-1.8,-1.8) {$c$};
\node (d2) at (1.8,-1.8) {$d$};

\draw (up) -- (a2);
\draw (up) -- (b2);
\draw (um) -- (c2);
\draw (um) -- (d2);

\draw (a2) -- (b2);
\draw (c2) -- (d2);

\draw[thick] (up) -- node[right, draw=none, shape=rectangle] {$e$} (um);

\node[draw=none, shape=rectangle] at (0,-2.5) {$H_u(\Sigma)$};
\end{scope}

\end{tikzpicture}
\caption{An example of the construction of $H_u(\Sigma)$. Solid edges in $\Sigma$ are positive, and dashed edges are negative.}\label{fig:Hu}
\end{figure}

Let $\Sigma$ be a signed graph with $n$ vertices.
Suppose that all negative edges of $\Sigma$ are incident to $u$.
Construct an ordinary graph $H_u(\Sigma)$ as follows:
replace $u$ by two vertices $u_+$ and $u_-$; replace each positive edge $uw$ by $u_+w$ and each negative edge $uw$ by $u_-w$;
keep every edge not incident to $u$; and add the edge $e=u_+u_-$. See Figure~\ref{fig:Hu}.

We call $e$ the distinguished edge of this construction.
The graph $H_u(\Sigma)$ is simple, even when $\Sigma$ has opposite signed parallel edges, and has $n+1$ vertices.

This is the usual splitting construction for a balancing vertex with an added edge between the split vertices.
For our purposes, we prove the following identity by counting colorings.
\begin{proposition}\label{prop:splitting}
With the notation above, we have
\begin{equation}\label{eq:splitting}
 P_{H_u(\Sigma)}(q)=(q-1)Z_\Sigma(q).
\end{equation}
Equivalently,
$\chi_{H_u(\Sigma)}(q)=(q-2)\chi_\Sigma(q)$.
\end{proposition}

\begin{proof}
Fix $a\in\calC_k$. From a proper zero-free coloring $c$ of $\Sigma$
with $c(u)=a$, define a coloring $\widetilde c$ of $H_u(\Sigma)$ by
\[
 \widetilde c(u_+)=a,\qquad
 \widetilde c(u_-)=-a,\qquad
 \widetilde c(w)=c(w)\quad(w\ne u).
\]
For a positive edge $uw$, its condition is exactly the proper coloring condition on $u_+w$;
for a negative edge, it is the condition on $u_-w$.
All remaining edges are positive and $\widetilde c(u_+) \neq \widetilde c(u_-)$.
Hence, this construction gives a bijection between the proper zero-free colorings of $\Sigma$ with $c(u)=a$ and the proper colorings of $H_u(\Sigma)$ by $\calC_k$ with $\widetilde c(u_+)=a$ and $\widetilde c(u_-)=-a$.

As in Lemma~\ref{lem:vertex-sum}(a), the number of colorings in the first set is $Z_\Sigma(2k)/(2k)$.
For $H_u(\Sigma)$, all ordered pairs of distinct colors on $u_+u_-$ occur equally often, so the number of colorings in the second set is $P_{H_u(\Sigma)}(2k)/(2k(2k-1))$.
Thus, \eqref{eq:splitting} holds for every $q=2k$, and hence as a polynomial identity.
\end{proof}

\begin{lemma}\label{lem:split-block}
Let $B$ be an unbalanced block whose negative edges are incident to a common vertex $u$. Then $H_u(B)$ is $2$-connected.
\end{lemma}
\begin{proof}
Since $B$ is $2$-connected, $B\setminus u$ is connected. Furthermore, $u$ is incident to both a positive edge and a negative edge. Indeed, if every edge incident to $u$ were negative, then every cycle of $B$ would contain an even number of negative edges, contradicting the unbalancedness of $B$.

Both $H_u(B) \setminus u_+$ and $H_u(B)\setminus u_-$ are connected,
since $B \setminus u$ is connected and each of $u_+$ and $u_-$ has a neighbor in $B\setminus u$.
For any $x\in V(B)\setminus\{u\}$, contracting the edge $u_+u_-$ in $H_u(B)\setminus x$ yields the underlying multigraph of $B\setminus x$, which is connected since $B$ is a block.
Thus, $H_u(B)\setminus x$ is also connected. Since $H_u(B)$ has at least three vertices, it is $2$-connected.
\end{proof}

For a connected ordinary graph $H$ with $N\ge2$ vertices, set \[ \beta(H)=(-1)^N P_H'(1). \]
This is known as the \emph{beta invariant} of its graphic matroid. It follows from \cite{Crapo} that
\begin{equation}\label{eq:beta-positive}
 \beta(H)\in\ZZ_{>0} \quad\text{if $H$ is $2$-connected or $H=K_2$.}
\end{equation}

\begin{corollary}\label{cor:block-evaluation}
Let $B$ be an unbalanced block with $m$ vertices whose negative edges are
incident to a common vertex $u$. Then \[ Z_B(1)=(-1)^{m+1}\beta(H_u(B)) \quad\text{and}\quad \beta(H_u(B))\ge1. \]
\end{corollary}
\begin{proof}
Differentiate \eqref{eq:splitting} and evaluate at $q=1$. Then we obtain $Z_B(1)=P_{H_u(B)}'(1)$.
Since $H_u(B)$ has $m+1$ vertices, the conclusion follows from Lemma~\ref{lem:split-block} and \eqref{eq:beta-positive}.
\end{proof}

\bigskip

\section{Proof of the main theorem}\label{sec:proof}

The goal of this section is to prove Theorem~\ref{thm:main}. 
\begin{proposition}\label{prop:negative}
Let $\Sigma$ be a connected signed graph on $n\ge2$ vertices whose negative edges are incident to a common vertex $u$.
Suppose that $\Sigma$ has no balanced block. Then
\[ Z_\Sigma(1)=(-1)^{n+1}b_\Sigma, \;\;\text{where }\; b_\Sigma=\prod_{B: \text{ a block of }\Sigma}\beta(H_u(B)) \in\ZZ_{>0}. \]
Consequently, the coefficients in \eqref{eq:first-coefficients} satisfy
\[
 \begin{array}{lll}
 n\text{ even}:& h_2\le-1,& h_3\ge n+1,\\[2pt]
 n\text{ odd}:& h_2\ge1,& h_3\le1-n.
 \end{array}
\]
In particular, $h_\Sigma(t)$ has a negative coefficient.
\end{proposition}
\begin{proof}
Since every block is unbalanced, it contains a negative edge and hence contains $u$.
Let $B_1,\ldots,B_b$ be the blocks with $m_i=|V(B_i)|$.
Evaluating \eqref{eq:block-product} at $q=1$ and using Corollary~\ref{cor:block-evaluation}, we obtain
\[ Z_\Sigma(1) =(-1)^{\sum_i(m_i+1)}   \prod_{i=1}^b\beta(H_u(B_i)).\]
By a decomposition into blocks, we see that $\sum_i(m_i-1)=n-1$, and hence $\sum_i(m_i+1)=n-1+2b$, which proves the desired identity. 
The coefficient bounds follow from \eqref{eq:first-coefficients} and $0\le Z_\Sigma(2)\le2$ by Lemma~\ref{lem:two-colors}.
For example, for $n$ even, we see that $h_2=Z_\Sigma(1)=(-1)^{n+1}b_\Sigma =-b_\Sigma \leq -1$ and
\[h_3=Z_\Sigma(2)-(n+1)Z_\Sigma(1) \geq 0+(n+1)b_\Sigma \geq n+1.\]
The case of odd $n$ is similar.
\end{proof}

\begin{proposition}\label{prop:realization}
Let $\Sigma$ be a connected signed graph with at least two vertices whose negative edges are incident to a common vertex.
If $\Sigma$ has a balanced block, then there is a connected ordinary graph $G$ on $|V(\Sigma)|$ vertices such that $P_G(q)=Z_\Sigma(q)$.
\end{proposition}
\begin{proof}
If every block is balanced, then every cycle is balanced, and the conclusion follows from \eqref{eq:balanced-ordinary}.
Otherwise, we may list the unbalanced blocks as $U_1,\ldots,U_r$ with $r\ge1$, and the balanced blocks as $B_1,\ldots,B_\ell$ with $\ell\ge1$.
Every $U_i$ contains the common endpoint $u$ of the negative edges. Distinct blocks have at most one common vertex,
so their union $\Sigma_0=U_1\cup\cdots\cup U_r$ is their vertex sum at $u$.
Let $H=H_u(\Sigma_0)$ with distinguished edge $e$. By Proposition~\ref{prop:splitting}, we have $P_H(q)=(q-1)Z_{\Sigma_0}(q)$.
For each balanced block $B_j$, let $G_j$ be a copy of its underlying ordinary graph, chosen so that $H,G_1,\ldots,G_\ell$ are pairwise vertex-disjoint. Then $P_{G_j}(q)=Z_{B_j}(q)$.
Choose any edge $f$ of $G_1$ and identify $f$ with $e$, including their endpoints, retaining one copy of the common edge.
Let $G_0$ be the resulting ordinary graph. Then it follows from \eqref{eq:edge-sum} that
\[ P_{G_0}(q) =\frac{P_H(q)P_{G_1}(q)}{q(q-1)} = \frac{Z_{\Sigma_0}(q)Z_{B_1}(q)}{q}.\]
Now attach $G_2,\ldots,G_\ell$ to $G_0$ by identifying one chosen vertex of each with a vertex of $G_0$.
The resulting graph $G$ is connected and simple, and we obtain the following by applying Lemma~\ref{lem:vertex-sum} repeatedly:
\[ P_G(q) =\frac{Z_{\Sigma_0}(q)\prod_{j=1}^\ell Z_{B_j}(q)}{q^\ell} =Z_\Sigma(q).\]
Thus, the balanced blocks can be reattached at arbitrary vertices; no assumption on their original positions is needed.

Finally, we have
\begin{align*}
 |V(G)| &=|V(H)|+|V(B_1)|-2 +\sum_{j=2}^\ell(|V(B_j)|-1)\\
 &=1+\sum_{i=1}^r(|V(U_i)|-1) +\sum_{j=1}^\ell (|V(B_j)|-1) =|V(\Sigma)|,
\end{align*}
as required.
\end{proof}

We are now ready to prove Theorem~\ref{thm:main}.

\begin{proof}[Proof of Theorem~\ref{thm:main}]
Since switching preserves all four conditions, we may assume that every negative edge of $\Sigma$ is incident to a common vertex.

\ref{main:graphic} $\Rightarrow$ \ref{main:positive}: 
Suppose that $Z_\Sigma(q)=P_G(q)$ for an ordinary graph $G$. Comparing degrees gives $|V(G)|=n$. Consequently, $h_\Sigma(t)=tW_G(t)$, and the assertion follows from \cite[Theorem~6]{Steingrimsson}.

\noindent
\ref{main:positive} $\Rightarrow$ \ref{main:block}:
Suppose that $\Sigma$ has no balanced block. By Proposition~\ref{prop:negative}, $h_\Sigma(t)$ has a negative coefficient, contradicting \ref{main:positive}. Thus, $\Sigma$ has a balanced block.

\noindent
\ref{main:block} $\Rightarrow$ \ref{main:graphic}:
This implication and the additional assertion about $G$ follow directly from Proposition~\ref{prop:realization}.

\noindent
\ref{main:block} $\Rightarrow$ \ref{main:uniform}:
Let $B$ be a balanced block of $\Sigma$. Since the underlying ordinary graph $\overline B$ has at least one edge, 
it follows from \eqref{eq:balanced-ordinary} that $Z_B(1)=P_{\overline B}(1)=0$. Hence, \eqref{eq:block-product} implies $Z_\Sigma(1)=0$. 

\noindent
\ref{main:uniform} $\Rightarrow$ \ref{main:block}:
If $\Sigma$ has no balanced block, then Proposition~\ref{prop:negative} says $Z_\Sigma(1)\ne 0$. 
\end{proof}



\bigskip

\section{Examples and the assumptions}\label{sec:examples}

In this section, we provide examples of signed graphs and their chromatic polynomials.
Some of them show that the assumptions of Theorem~\ref{thm:main} cannot be omitted.

\begin{example}[Signed cycles]\label{ex:cycles}
Let $\Sigma$ be a signed cycle $C_n$ with $n\ge3$.
If the number of negative edges is even, then $\Sigma$ is balanced.
If it is odd, then $\Sigma$ can be switched to a signed cycle with exactly one negative edge.
By splitting one of its endpoints, we obtain $C_{n+1}$. Therefore, the formula for ordinary cycles and \eqref{eq:splitting} provide 
\[ Z_\Sigma(q)= \begin{cases}
 (q-1)^n+(-1)^n(q-1),&\text{if $\Sigma$ is balanced},\\
 (q-1)^n-(-1)^n,&\text{if $\Sigma$ is unbalanced}.
 \end{cases}\]
Equivalently,
\[ \chi_\Sigma(q)=
 \begin{cases}
 (q-2)^n+(-1)^n(q-2),&\text{if $\Sigma$ is balanced},\\
 (q-2)^n-(-1)^n,&\text{if $\Sigma$ is unbalanced}.
 \end{cases}
\]
For an unbalanced $C_3$, this gives 
\[ \chi_\Sigma(q)=(q-2)^3+1 \:\;\text{and}\;\; h_\Sigma(t)=t^2-2t^3+7t^4.\] 
For a balanced $C_3$, we have $Z_\Sigma(q)=q(q-1)(q-2)$ and $h_\Sigma(t)=6t^4$. 
Although $\Sigma$ is 2-connected, all four conditions of Theorem~\ref{thm:main} are satisfied. This illustrates that Theorem~\ref{thm:main} applies to both balanced and unbalanced signed graphs.
\end{example}

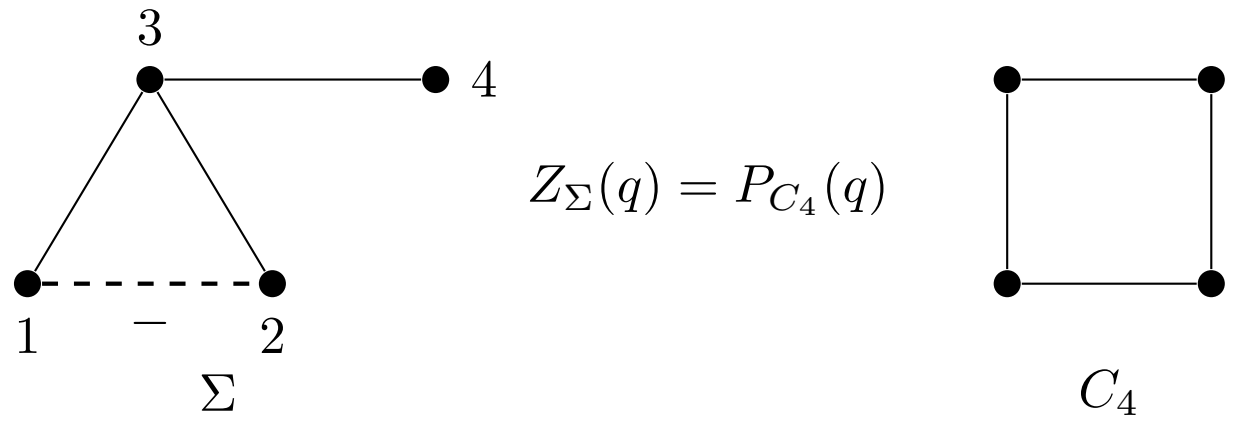
\begin{figure}[htbp]
\centering
\begin{tikzpicture}[scale=.9,every node/.style={font=\small},
 vertex/.style={circle,fill,inner sep=1.8pt}]
 \begin{scope}
  \node[vertex,label=below:$1$] (a) at (0,0) {};
  \node[vertex,label=below:$2$] (b) at (1.8,0) {};
  \node[vertex,label=above:$3$] (c) at (.9,1.5) {};
  \node[vertex,label=right:$4$] (d) at (3,1.5) {};
  \draw[dashed,thick] (a)--node[below]{$-$}(b);
  \draw (a)--(c)--(b) (c)--(d);
  \node at (1.4,-.8) {$\Sigma$};
 \end{scope}
 \node at (5,0.7) {$Z_\Sigma(q)=P_{C_4}(q)$};
 \begin{scope}[xshift=7.2cm]
  \node[vertex] (v1) at (0,0) {};
  \node[vertex] (v2) at (1.5,0) {};
  \node[vertex] (v3) at (1.5,1.5) {};
  \node[vertex] (v4) at (0,1.5) {};
  \draw (v1)--(v2)--(v3)--(v4)--(v1);
  \node at (.75,-.8) {$C_4$};
 \end{scope}
\end{tikzpicture}
\caption{An unbalanced $C_3$ with a balanced bridge, and an ordinary graph with the same zero-free chromatic polynomial. The dashed edge is the negative edge.}\label{fig:bridge}
\end{figure}

\begin{example}[A balanced bridge]\label{ex:bridge}
Attach a pendant vertex by a positive edge to an unbalanced $C_3$, and let $\Sigma$ be the resulting signed graph.
The pendant edge is a balanced block, and
\[ Z_\Sigma(q) =(q-1)\bigl((q-1)^3+1\bigr) =(q-1)^4+(q-1) =P_{C_4}(q).\]
Accordingly,
\[ \chi_\Sigma(q)=P_{C_4}(q-1), \qquad h_\Sigma(t)=2t^3+8t^4+14t^5.\]
This is the simplest illustration of how a balanced block removes the extra factor in the splitting identity; see Figure~\ref{fig:bridge}.
\end{example}

\begin{example}[No balancing vertex]\label{ex:no-balancing-vertex}
Let $\Sigma$ have the vertex set $\{1,2,3,4,5\}$ and
\[ E^+=\{13,14,23,24,35\}, \qquad E^-=\{12,34\}.\]
Its underlying graph is $K_4$ with a pendant edge; see Figure~\ref{fig:counterexample}.
In particular, it is connected and not $2$-connected, and the edge $35$ is a balanced block.
All four triangles on $\{1,2,3,4\}$ are unbalanced. Therefore, deleting any vertex leaves a graph containing an unbalanced triangle, 
i.e., $\Sigma$ has no balancing vertex, even after switching.
\begin{figure}[htbp]
\centering
\begin{tikzpicture}[scale=0.9,every node/.style={font=\small},
 vertex/.style={circle,fill,inner sep=1.8pt}]
 \node[vertex,label=below left:$1$] (a) at (0,0) {};
 \node[vertex,label=below right:$2$] (b) at (1.8,0) {};
 \node[vertex,label=above:$3$] (c) at (1.8,1.8) {};
 \node[vertex,label=above left:$4$] (d) at (0,1.8) {};
 \node[vertex,label=right:$5$] (e) at (3.4,1.8) {};
 \draw (a)--(d) (a)--(c) (b)--(c) (b)--(d) (c)--(e);
 \draw[dashed,thick] (a)--node[below]{$-$}(b);
 \draw[dashed,thick] (d)--node[above]{$-$}(c);
\end{tikzpicture}
\caption{A connected signed graph with a balanced block and a
nonnegative numerator, whose zero-free chromatic polynomial is not
an ordinary chromatic polynomial. Negative edges are dashed.}
\label{fig:counterexample}
\end{figure}
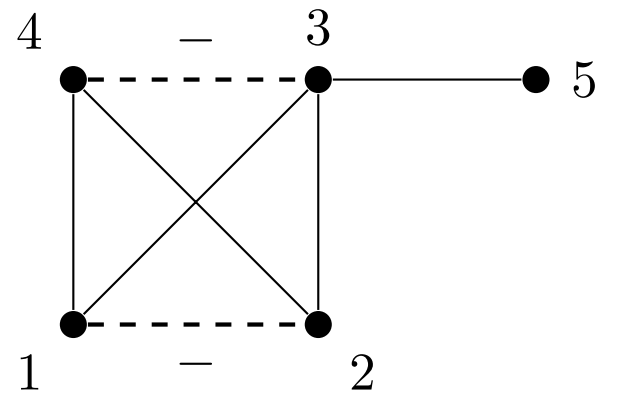
Switching at vertices $1$ and $2$ makes all six edges of $K_4$ negative, while $35$ is still positive.
Let $K_4^-$ be $K_4$ having all negative edges and partition its four vertices into nonempty color classes. 
For a fixed partition into $j$ classes, the number of ways to assign distinct colors with no opposite pair is $\prod_{i=0}^{j-1}(q-2i)$, for even $q$.
The numbers of partitions into $j=1,2,3,4$ classes are $1,7,6,1$, respectively. It follows that
\begin{align*}
 Z_{K_4^-}(q) &=q+7q(q-2)+6q(q-2)(q-4)+q(q-2)(q-4)(q-6)\\
 &=q^4-6q^3+15q^2-13q.
\end{align*}
The pendant vertex contributes a factor $q-1$, giving the following:
\begin{equation}\label{eq:counterexample-Z}
 Z_\Sigma(q)=q^5-7q^4+21q^3-28q^2+13q.
\end{equation}
Consequently, we obtain
\begin{align*}
 \chi_\Sigma(q) &=q^5-12q^4+59q^3-143q^2+165q-70, \\
 h_\Sigma(t) &=2t^3+18t^4+30t^5+70t^6.
\end{align*}
In particular, the numerator is coefficientwise nonnegative.

Nevertheless, $Z_\Sigma$ is not the chromatic polynomial of an ordinary graph. Suppose otherwise and let $G$ be such a graph.
Its first two coefficients force $G$ to have five vertices and seven edges. By inclusion--exclusion,
the coefficient of $q^3$ in $P_G(q)$ is $\binom72-\tau$, where $\tau$ is the number of triangles in $G$.
Therefore, \eqref{eq:counterexample-Z} would imply $\tau=\binom72-21=0$.
However, a triangle-free graph on five vertices has at most six edges, a contradiction.
Thus, this example satisfies conditions \ref{main:positive}, \ref{main:block}, and \ref{main:uniform} of Theorem~\ref{thm:main}, but not \ref{main:graphic}. 
This shows that the balancing-vertex assumption cannot be omitted.
\end{example}

\begin{example}[Connectedness]\label{ex:disconnected}
Let $\Sigma$ be the graph consisting of an unbalanced $C_3$ with one negative edge and two isolated vertices.
Then \[ Z_\Sigma(q)=q^3(q^2-3q+3), \qquad \chi_\Sigma(q)=(q-1)^3(q^2-5q+7),\] and
\[ h_\Sigma(t)=t^2+2t^3+48t^4+62t^5+7t^6\]
has nonnegative coefficients.
On the other hand, $Z_\Sigma(1)=1$. Any ordinary graph with an edge has $P_G(1)=0$, while an edgeless graph on five vertices has polynomial $q^5$, different from $Z_\Sigma(q)$.
Therefore, no ordinary graph has chromatic polynomial $Z_\Sigma(q)$.
This shows that connectedness cannot be omitted even if only the equivalence between \ref{main:graphic} and \ref{main:positive} is retained.
\end{example}

\bigskip

\section*{Acknowledgements}
The first named author was supported by JSPS KAKENHI Grant Numbers JP24K00521 and JP23H00081.

\end{document}